\documentclass[11pt]{article}

\usepackage{amsmath}
\usepackage{amsthm}
\usepackage{lineno}
\usepackage{graphicx}
\usepackage{caption}
\usepackage{subcaption}
\usepackage[a4paper,total={6in,8in}]{geometry}
\title{Specular reflection on the surface of a sphere: compass and ruler constructions.}
\author{Nikolaos K. Kollas\footnote{kollas@upatras.gr}\\\small{Division of Theoretical and Mathematical Physics}\\ \small{Physics Department} \small{University of Patras, Greece}}
\date{}
\begin{document}
\maketitle
%
\abstract{We provide an explicit \emph{geometric algorithm} involving only ruler and compass constructions in order to specify the \emph{specular reflection} point on the surface of a reflecting sphere of radius $r$ given two \emph{focal points} $A$ and $B$ lying outside of it. By numerically implementing the algorithm we compute the point in question for a number of cases. We conclude by discussing how the first iteration of the algorithm constitutes a \emph{first order approximation} to the real solution by providing a closed expression for it as well as the error involved in doing so, as a function of the distances of the two focal points from the origin and the angle formed between them.}
\section{Introduction.}
A point of \emph{specular reflection} is defined as the unique point on a reflective surface such that the angles formed by the \emph{incident} and \emph{reflecting rays} from two \emph{focal points} lying outside the surface with the perpedincular direction at that point become equal.

The problem of specifying the specular point in the case where the surface is that of a sphere appeared for the first time in \cite{Fishback}. This is of particular interest in the field of \emph{terrestrial} and \emph{satellite communications}, where knowledge of the exact point is necessary in order to avoid destructive interference between the direct signal sent over the air from \emph{transmitter} to \emph{receiver} and the one reflected from the surface of the earth. Although exact analytic solutions are known to exist \cite{Allen}, their complexity dramatically reduces our ability of obtaining a clear intuitive picture of the problem in contrast with a solution which would be based on geometric arguments. 

In this letter we introduce such a solution by providing a \emph{geometric algorithm} based on \emph{ruler} and \emph{compass} constructions. This approach constitutes a marked improvement to the analytic method especially when a simple answer is needed in order to make quick qualitative descriptions in applications, in this case the first iteration of the algorithm provides a satisfying approximation.

In Section 2 we briefly state the problem of specular reflection on the surface of a sphere and discuss how to construct the relevant equations. The main part consists of Section 3 where we explicitly describe the algorithm and provide a proof for it's correctness as well as numerical results for a number of cases. Section 4 is concerned with the error made in approximating the exact solution with the first iteration of the algorithm for which a closed analytic expression is shown to exist while Section 5 includes closing remarks. 
\section{Analytic construction.}
Let $A$ and $B$ be two focal points lying outside the surface of a sphere of radius $r$ centered at $O$ at distances $r_A$ and $r_B$ respectively from the origin forming an angle of $2\theta$ as in Figure 1 where without loss of generality we always assume that $r_A\leq r_B$. We are interested in specifying the point of \emph{specular reflection} $M$ lying on the surface of the sphere such that the angles formed between $MA$, $MB$ and the ray connecting the origin with point $M$ become equal. Since this point lies on the \emph{reflection plane} defined by the incident and reflected rays it suffices to restrict our attention to the circle of radius $r$ formed by the intersection of the sphere with this surface.

Any possible attempt of attacking the problem will inevitably require the solution of a \emph{quartic polynomial equation}. For example, it should be immediately obvious from the above figure that the two triangles $(AMC)$ and $(BMD)$ formed by taking the distances of $A$ and $B$ to the ray connecting the origin to point $M$ are similar, if $\phi$ is it's angular distance from the bisector of $\widehat{AOB}$ it follows that
\begin{figure}\label{fig1}
	\centering
		\includegraphics[scale=0.6]{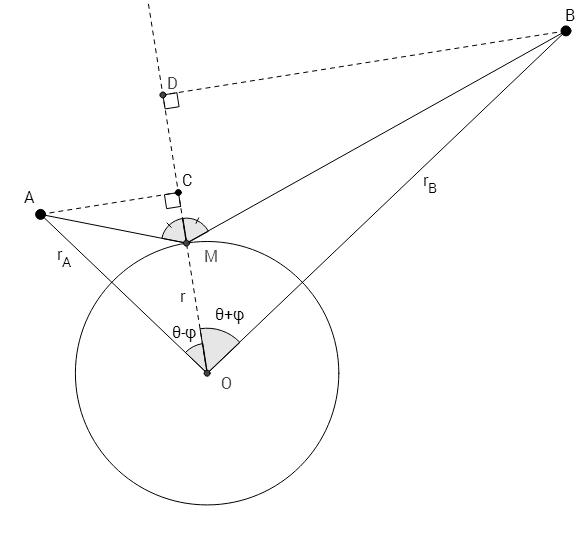}
		\caption{Geometry of specular reflection on the surface of a sphere given two focal points $A$ and $B$. (Note how the two angles at point $O$ sum to $2\theta$).}
\end{figure}
	$$\frac{r_A\sin(\theta-\phi)}{r_A\cos(\theta-\phi)-r}=\frac{r_B\sin(\theta+\phi)}{r_B\cos(\theta+\phi)-r}$$
which after a little bit of algebra can be rewritten as
\begin{equation}\label{eq1}
	(1+\tan^2\phi)\left[\left(\frac{r}{r_A}-\frac{r}{r_B}\right)\tan\theta+\left(\frac{r}{r_A}+\frac{r}{r_B}\right)\tan\phi\right]^2=4\tan^2\phi(1+\tan^2\theta)
\end{equation}
Equation (\ref{eq1}) is a quartic polynomial with respect to $\tan\phi$, although exact solutions can be found either with the use of a symbolic algebraic program, the general method for solving quartic equations developed by Cardano \cite{cardano} or by direct analytic methods (see for instance \cite{Allen} for the treatment of a similar quartic equation with complex coefficients), these tend to be extremely complicated (except for the case where $r_A=r_B$)\footnote{In this case the specular reflection point lies on the intersection of the circle of radius $r$ with the bisector of the angle between $A$ and $B$ $(\phi=0)$.} and require the use of numerical methods. Moreover it can be shown that in general all four solutions are real and distinct (see Appendix A for a proof and a physical interpretation). Since it is virtually impossible to tell which one corresponds to point $M$ we are forced to calculate all four of them and then choose the one lying on the arc of the circle contained in $\widehat{AOB}$. As this seems highly impractical, an alternative solution to the problem is therefore needed.

\section{Geometric construction.}
We shall now explicitly describe a \emph{geometric algorithm} involving only ruler and compass constructions for the specification of the \emph{specular reflection point}. The algorithm consists of the following steps
\begin{enumerate}
	\item Make an initial guess for the point of specular reflection on the sphere $M_0=(r,\phi_0)$.
	\item Connect this point to point $B$. Find the intersection $A'$ of $M_0B$ with the circle of radius $r_A$ centered at $O$.
	\item Replace $M_0$ with $M'=(r,\phi')$ where $\phi'$ is such that the point lies on the bisector of the angle $\widehat{AOA'}$.
	\item Repeat steps 2 and 3 $n$ times.
\end{enumerate}

Figure 2 presents the first iteration of the algorithm $(n=0)$.
\begin{figure}\label{fig2}
	\centering
		\begin{subfigure}{0.3\textwidth}
		\includegraphics[scale=0.25]{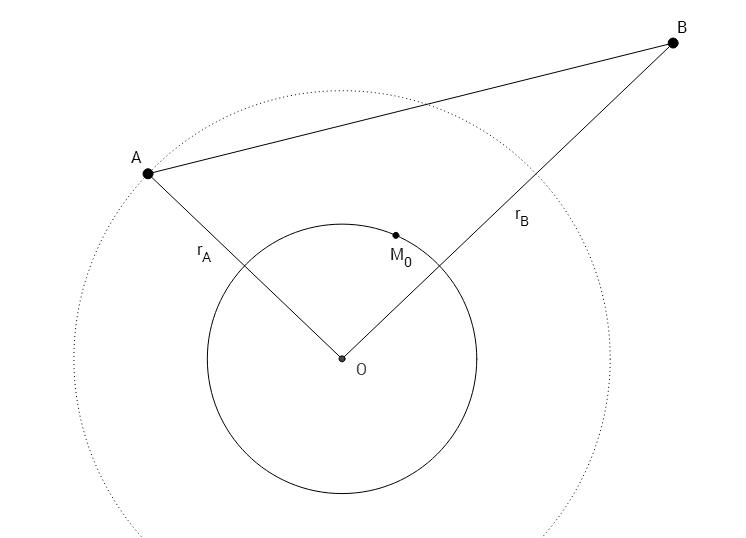}
		\caption{Step 1}
		\end{subfigure}
		\begin{subfigure}{0.3\textwidth}
		\includegraphics[scale=0.25]{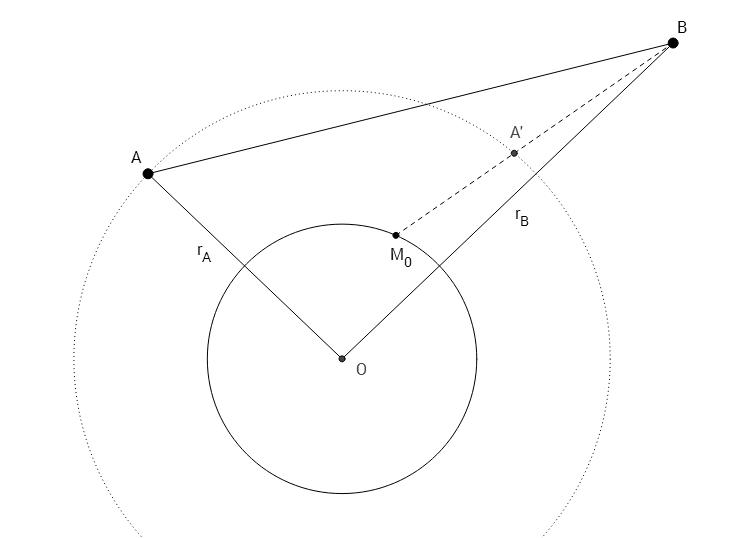}
		\caption{Step 2}
		\end{subfigure}
		\begin{subfigure}{0.3\textwidth}
		\includegraphics[scale=0.25]{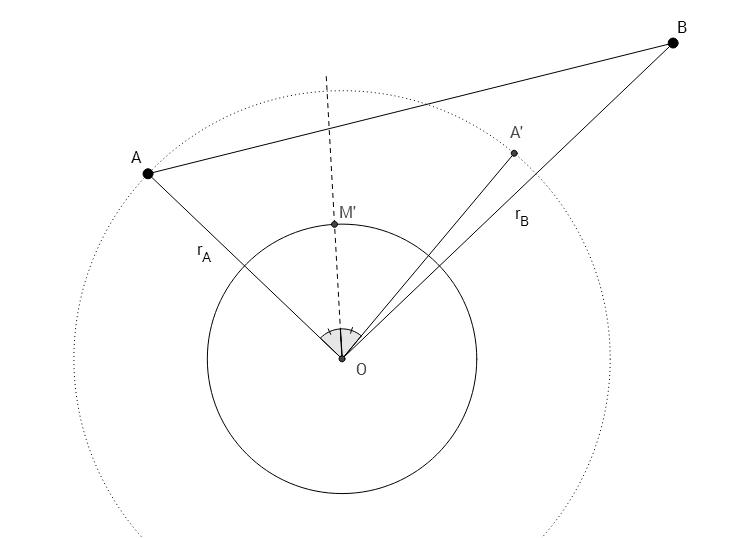}
		\caption{Step 3}
		\end{subfigure}
\caption{First iteration of the geometric algorithm.}
\end{figure}
We will now prove that in the limit of an infinite number of iterations $(n\to\infty)$, the sequence of points generated by the algorithm converges to $M$. 
\begin{proof}
Let $f(\phi)$ be the function that outputs the polar angle of the first iteration of the algorithm given an initial guess $\phi$. We note that the point of specular reflection is a fixed point of $f$, namely suppose $\phi_M$ is it's polar angle then $f(\phi_M)=\phi_M$. Furthermore Equation \ref{eq1} suggests that it must lie on the arc of the circle bounded between the bisector of $\widehat{AOB}$ and the perpendicular distance of $AB$ from the origin. If we can thus establish that $f$ is a \emph{contraction mapping} \footnote{A contraction mapping is a function from a metric space $X$ to itself such that for any two elements $x,y\in X$ there exists a nonnegative real number $k<1$ such that $d(f(x),f(y))\leq k d(x,y)$ with $d$ the metric on $X$.} on this interval then a simple application of the \emph{Banach fixed point theorem} \cite{debnath} suffices to prove that this is the unique fixed point. 

To show that $f$ is a \emph{contraction} suppose $d(x,y)=|x-y|$ is the metric on the unit circle with $x,y\in[-\pi,\pi]$ and let $\phi_0$ and $\psi_0$ be two initial guesses in it's domain with $\phi_0>\psi_0$ as in Figure 3, with the bisector represented by the dashed and dotted line, from which it can also be seen that $f$ is necessarily a monotonically increasing function. It is now true that  $f(\phi_0)-f(\psi_0)=\frac{\phi'-\psi'}{2}$, but because the distance between points $(r,\phi_0)$, $(r_A,\phi')$ is greater than that of $(r,\psi_0)$, $(r_A,\psi')$ it follows that $\psi_0-\psi'<\phi_0-\phi'$ so $f(\phi_0)-f(\psi_0)<\frac{\phi_0-\psi_0}{2}$ since the choice of $\phi_0$ and $\psi_0$ was arbitrary this must hold for any two choices in the interval.
\end{proof}
\begin{figure}
	\centering
		\includegraphics[scale=0.35,trim = 0mm 0mm 0mm 9cm]{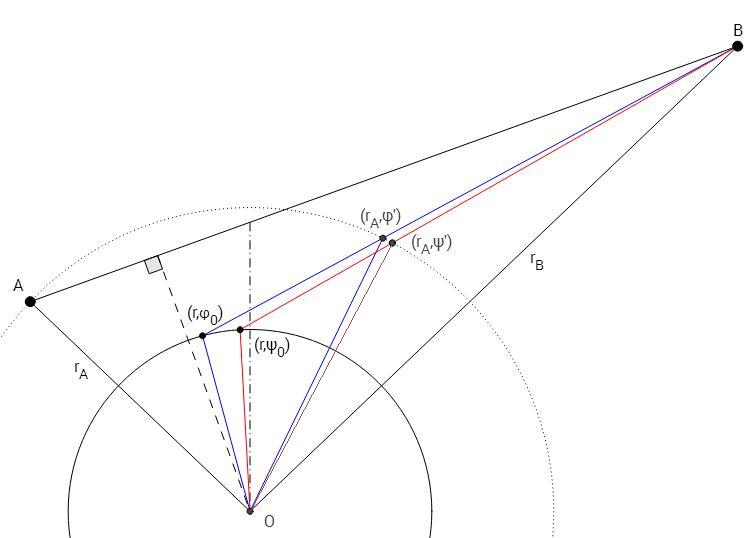}
	\caption{}
\end{figure}
In Figure 4 we plot the exact point of specular reflection as a function of $\theta$ for a number of distances of the two focal points, obtained from a numerical implementation of the geometric algorithm as a simple loop in a Python script. As was already mentioned in Section 2, angular distances are measured with respect to the bisector of $\widehat{AOB}$, this means that $0\leq\phi\leq\theta+\cos^{-1}\frac{r^2_{AB}+r^2_A-r^2_B}{2r_Ar_{AB}}-\frac{\pi}{2}$ where $r_{AB}$ is the distance between the two focal points.

Figures 5 (a)-(d) display the output of the first iteration $f(\phi)$ for fixed values of $r_A$, $r_B$ and $\theta$ from which the specular reflection point is easily determined by simply demanding that $f(\phi)=\phi$. Note that since for a fixed value of the distances, the angle formed between them cannot increase more than the value it will attain when the distance of point $A$ from $B$ ($r_{AB}$) is tangent to the sphere, $\theta$ will range between $0\leq2\theta\leq\cos^{-1}\frac{r^2-\sqrt{(r_A^2-r^2)(r_B^2-r^2)}}{r_Ar_B}$. Numerical results suggest that the sequence of approximations to the specular point \emph{converge linearly} to the real point.
\begin{figure}
	\centering
		\includegraphics[width=6.8in,trim = 25mm 80mm 0mm 14cm]{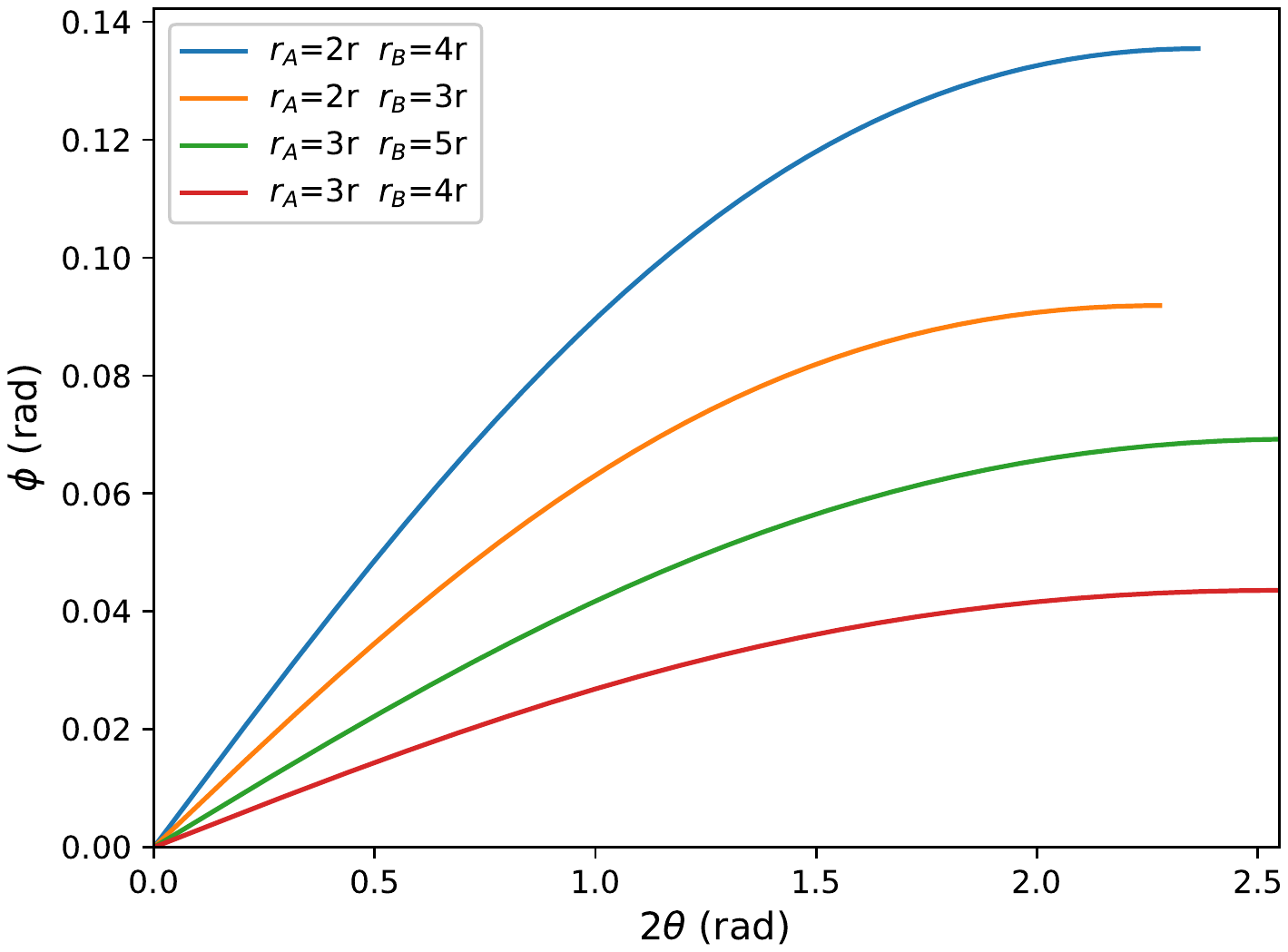}
		\caption{Point of specular reflection as a function of the angle $\widehat{AOB}$.}
		\begin{subfigure}{0.4\textwidth}
			\includegraphics[width=2.8in,trim = 25mm 80mm 0mm 90mm]{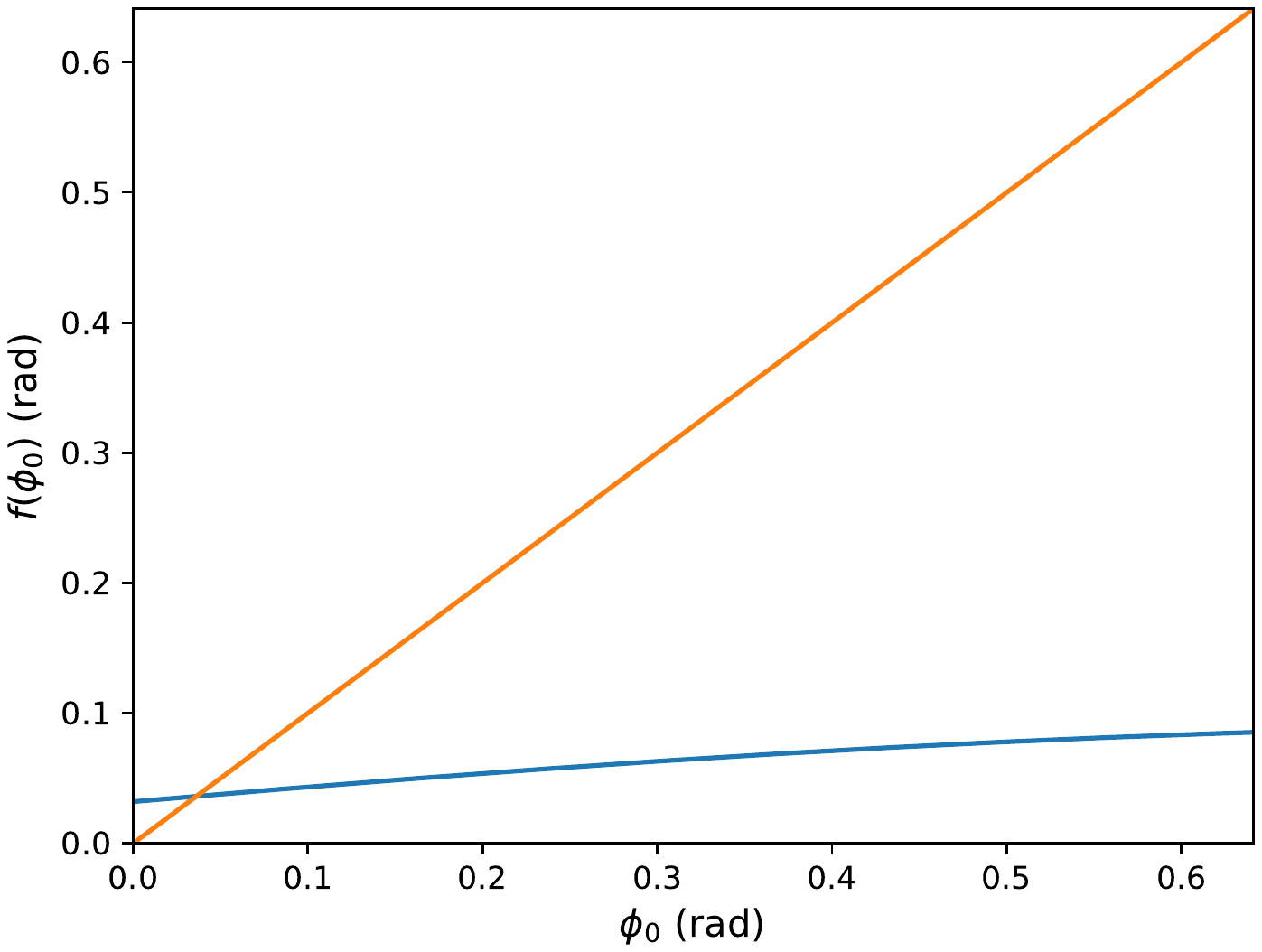}
			\caption{$r_A=2r$, $r_B=3r$, $2\theta=\frac{\pi}{6}$}
		\end{subfigure}
		\begin{subfigure}{0.4\textwidth}
			\includegraphics[width=2.8in,trim = 25mm 80mm 0mm 90mm]{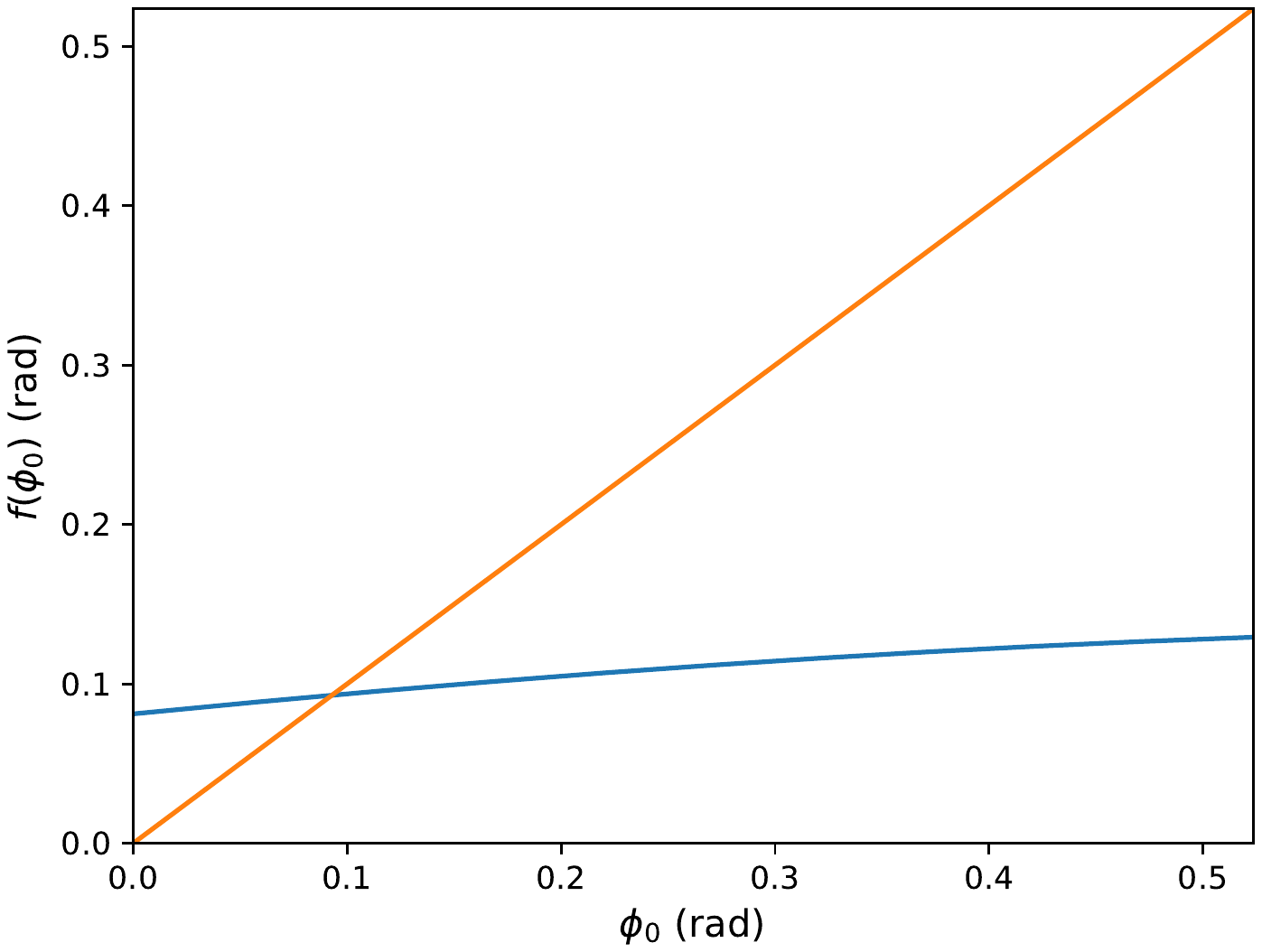}
			\caption{$r_A=2r$, $r_B=4r$, $2\theta=\frac{\pi}{3}$}
		\end{subfigure}
		\begin{subfigure}{0.4\textwidth}
			\includegraphics[width=2.8in,trim = 25mm 80mm 0mm 90mm]{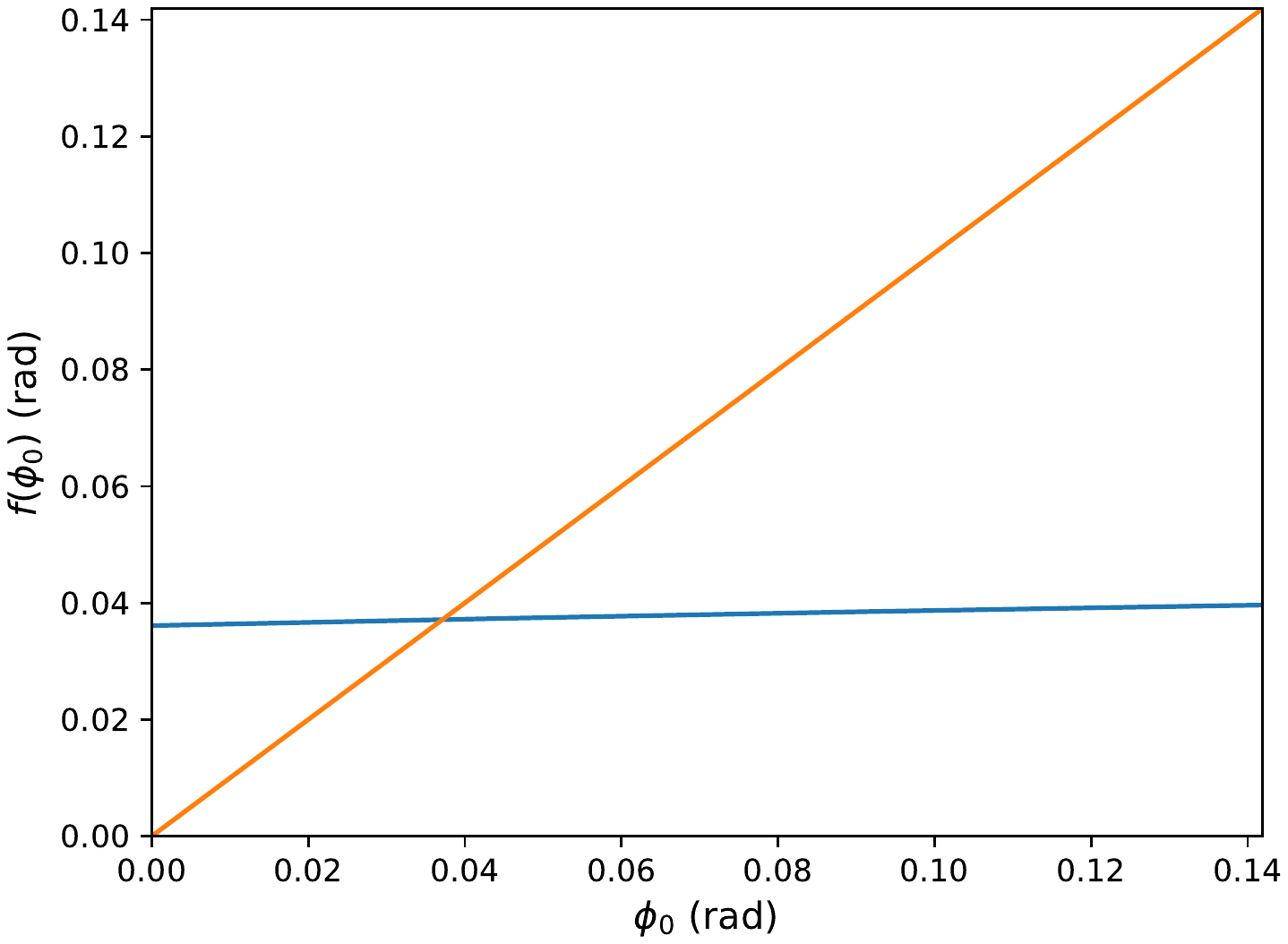}
			\caption{$r_A=3r$, $r_B=4r$, $2\theta=\frac{\pi}{2}$}
		\end{subfigure}
		\begin{subfigure}{0.4\textwidth}
			\includegraphics[width=2.8in,trim = 25mm 80mm 0mm 90mm]{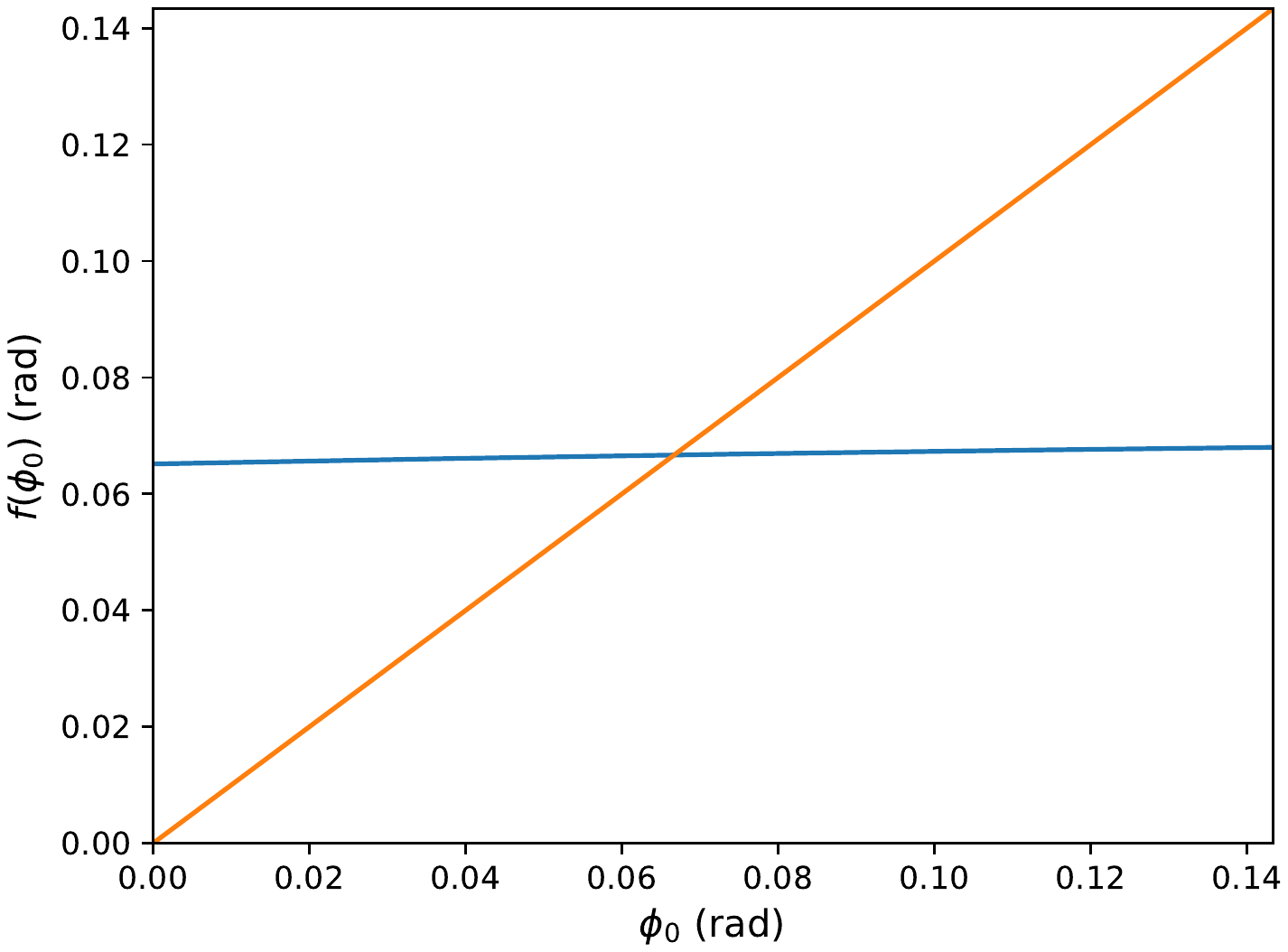}
			\caption{$r_A=3r$, $r_B=5r$, $2\theta=\frac{2\pi}{3}$}
		\end{subfigure}
	\caption{}
\end{figure}
\section{First order approximation.}
Let $M_0=(r,\phi_0)$ be our initial guess and $A'=(r_A,-\phi_{A'})$, $M'=(r,\frac{\theta-\phi_{A'}}{2})$ the points generated after the first iteration of the algorithm $(n=0)$ as in Figure 6 where as in Figure 3 the bisector is represented by the dashed and dotted line. Since the area $E_B$ of $(OM_0B)$ is equal to that of $(OM_0A')$ and $(OA'B)$ it follows that $rr_B\sin(\theta+\phi_0)=r_Ar_B\sin(\theta-\phi_{A'})+rr_Asin(\phi_0+\phi_{A'})$. Solving for $\theta-\phi_{A'}$ we find that
\begin{equation}\label{eq2}
	\frac{\theta-\phi_{A'}}{2}=\tan^{-1}\frac{2E_B\sqrt{r_A^2r_{M_0B}^2-4E_B^2}-2E_A\sqrt{r_A^2r_{M_0B}^2-4E_A^2}}{r_A^2r_{M_0B}^2-4E_B^2-4E_A^2}
\end{equation}
\begin{figure}
	\includegraphics[scale=0.45]{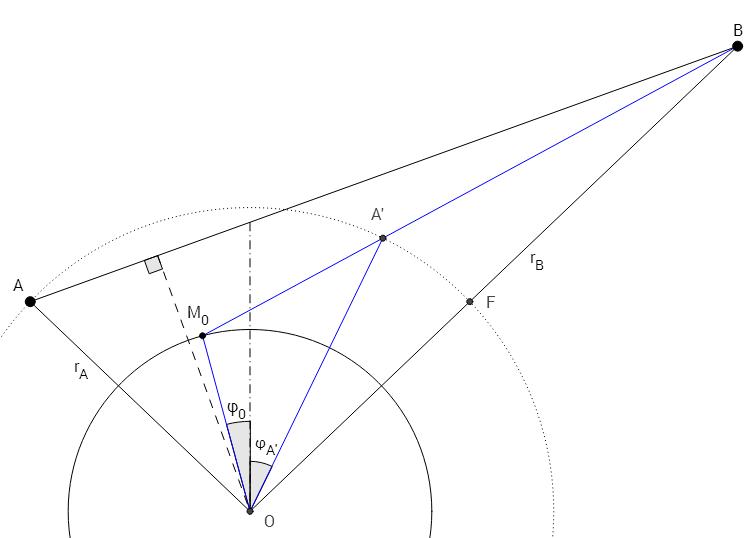}
	\caption{}
\end{figure}
where $E_A$ is the area of $(OM_0F)$, $F$ being the intersection of $OB$ with the circle of radius $r$ and $r_{M_0B}$ the distance of the initial guess from point $B$. Point $M'$ can be considered as a \emph{first order approximation} to the specular reflection point.

In Figure 7 we plot the \emph{relative error} to the first order approximation for the polar angle of the specular reflection point for an initial guess of $\phi_0=0$ as a function of the angle $\widehat{AOB}$ for varying distances of the two focal points. We note that in these cases the approximation is accurate within three orders of magnitude. Equation \ref{eq2} constitutes therefore a satisfying analytic estimation in cases where a quick and qualitative solution is needed.

In Figures 8 (a)-(d) we present the relative error as a function of the initial guess. As is expected this decreases the closer our initial guess is to the actual point of specular reflection.
\begin{figure}
	\centering
		\includegraphics[width=6.8in,trim = 25mm 80mm 0mm 14cm]{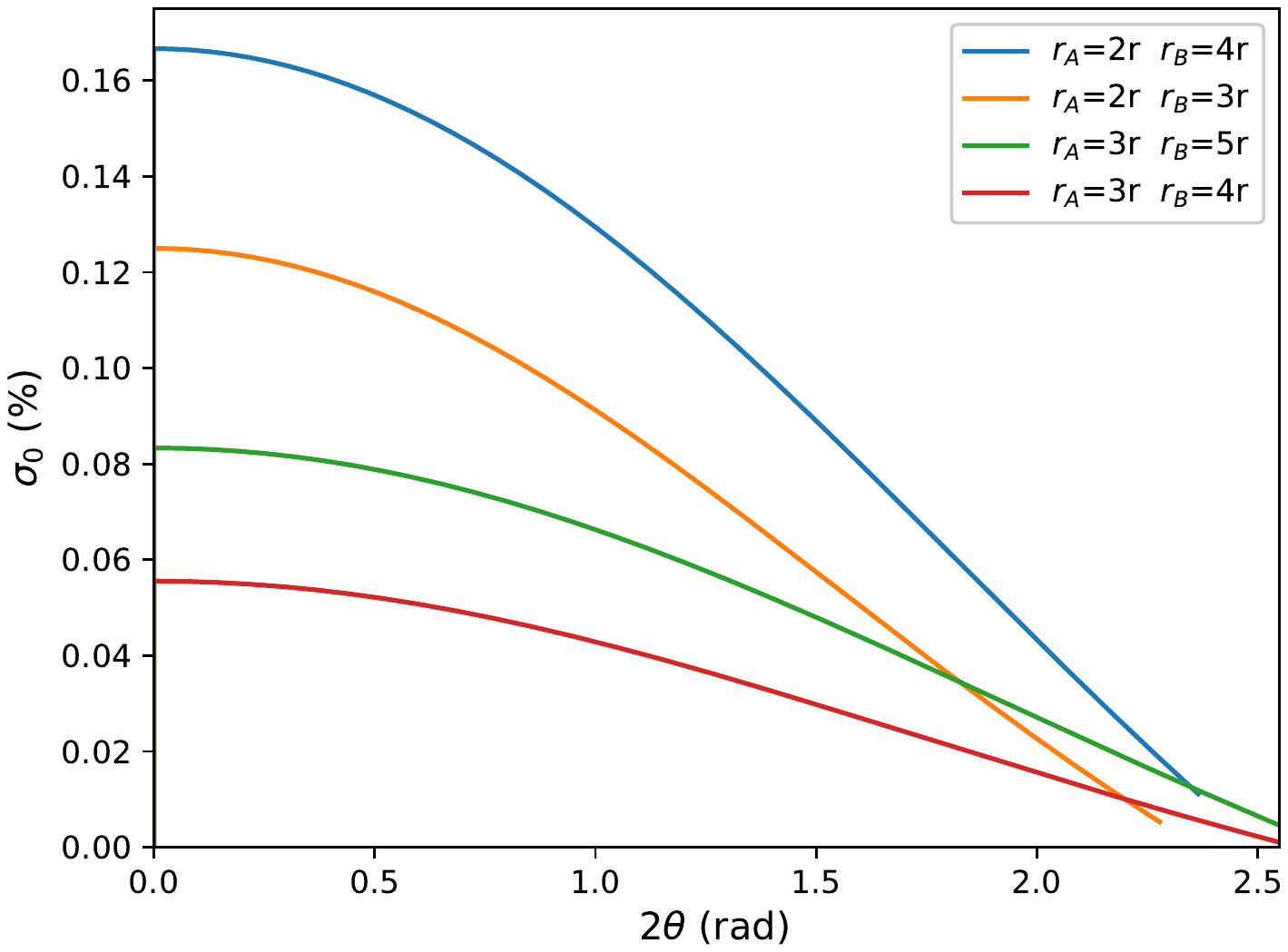}
		\caption{Relative error of the first order approximation as a function of the angle $\widehat{AOB}$.}
		\begin{subfigure}{0.4\textwidth}
			\includegraphics[width=2.8in,trim = 25mm 80mm 0mm 90mm]{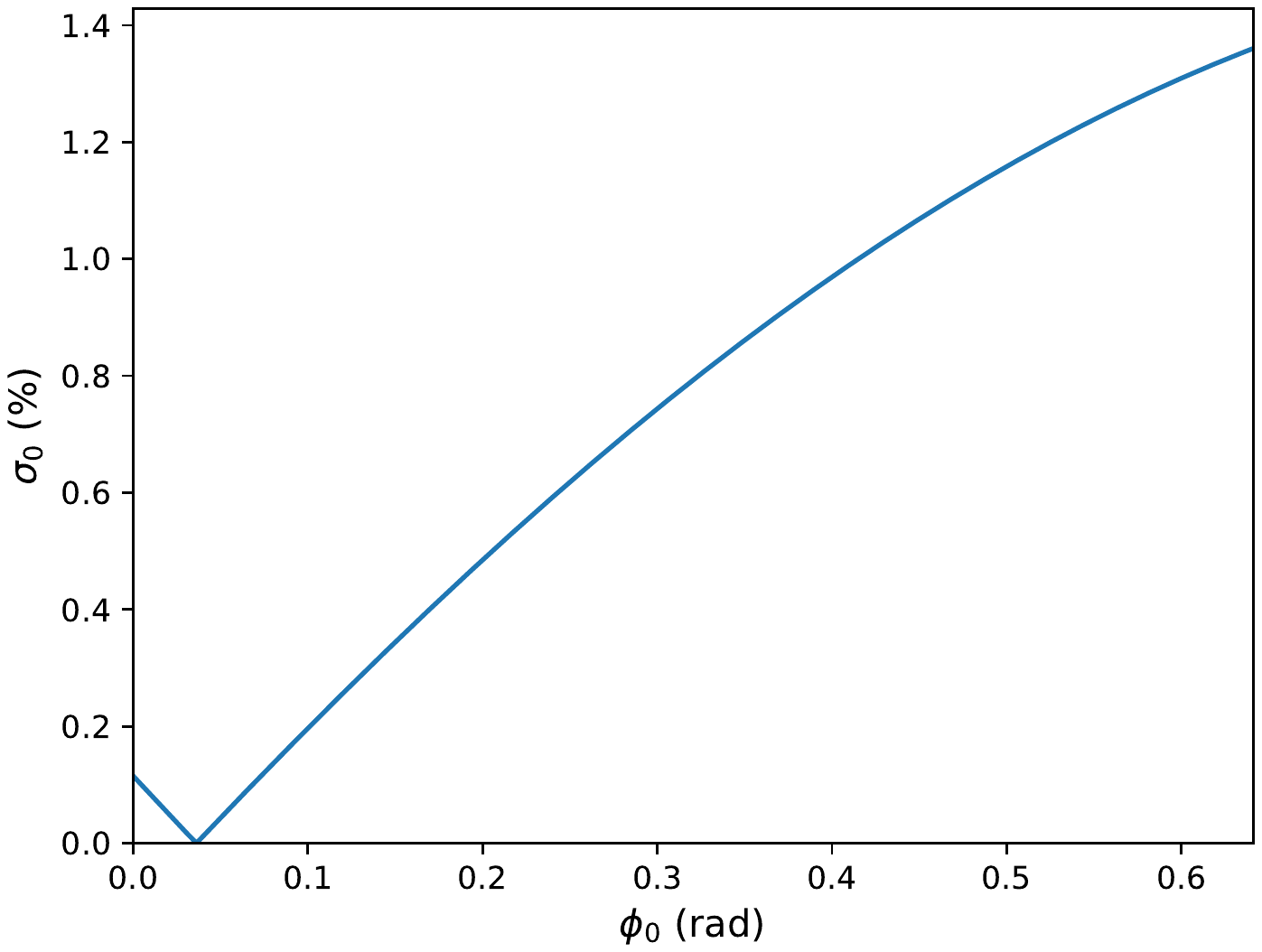}
			\caption{$r_A=2r$, $r_B=3r$, $2\theta=\frac{\pi}{6}$}
		\end{subfigure}
		\begin{subfigure}{0.4\textwidth}
			\includegraphics[width=2.8in,trim = 25mm 80mm 0mm 90mm]{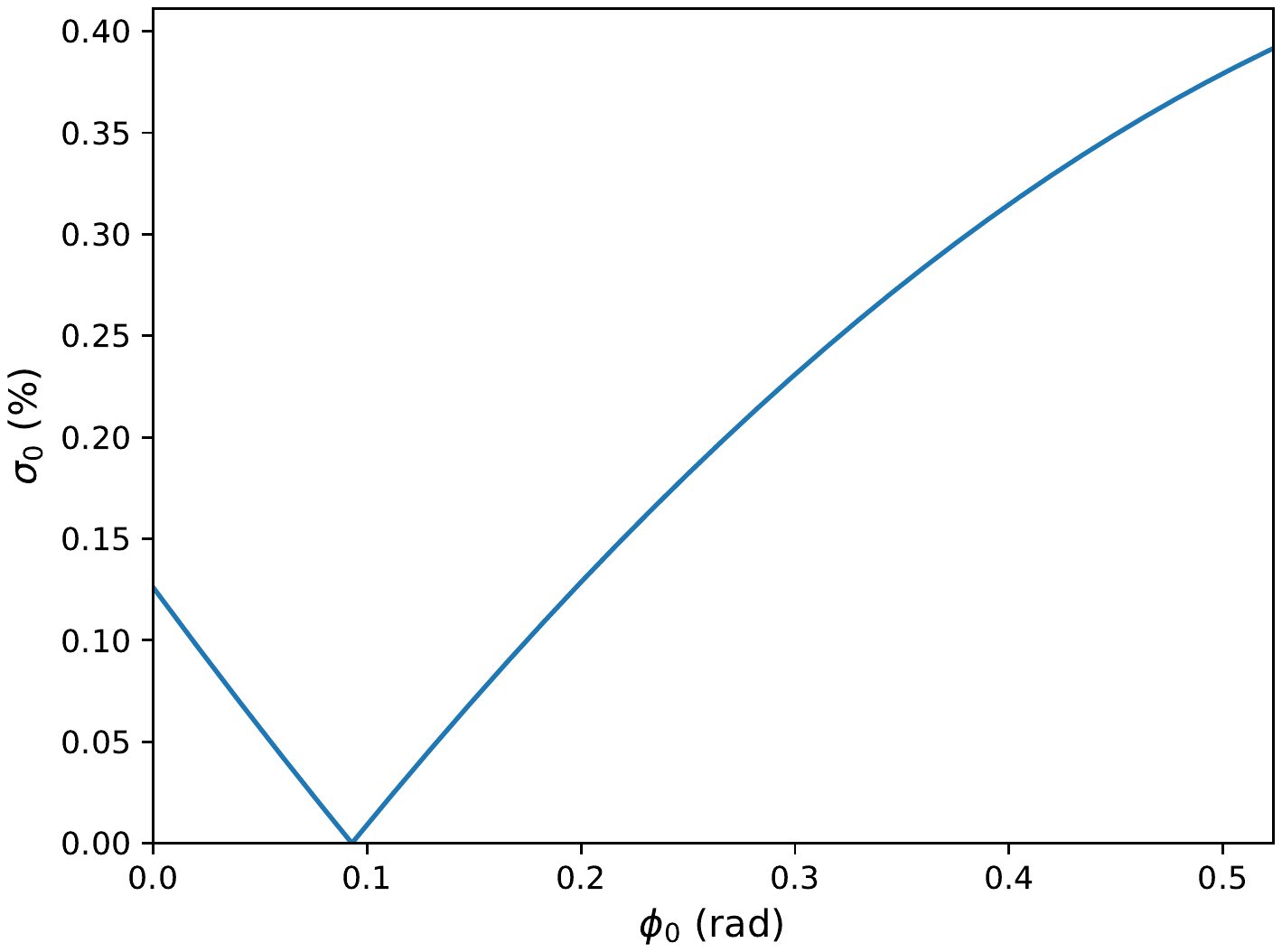}
			\caption{$r_A=2r$, $r_B=4r$, $2\theta=\frac{\pi}{3}$}
		\end{subfigure}
		\begin{subfigure}{0.4\textwidth}
			\includegraphics[width=2.8in,trim = 25mm 80mm 0mm 90mm]{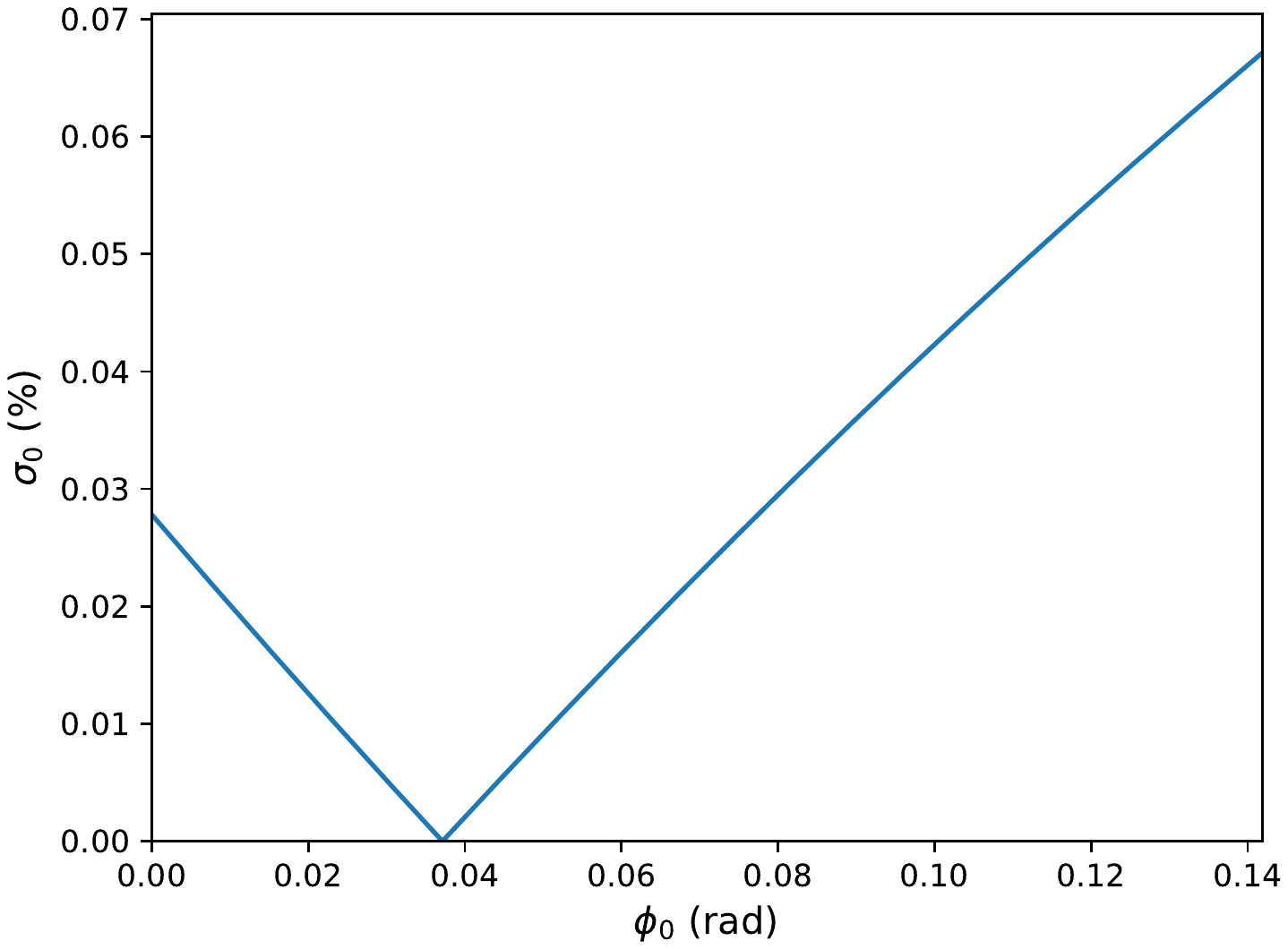}
			\caption{$r_A=3r$, $r_B=4r$, $2\theta=\frac{\pi}{2}$}
		\end{subfigure}
		\begin{subfigure}{0.4\textwidth}
			\includegraphics[width=2.8in,trim = 25mm 80mm 0mm 90mm]{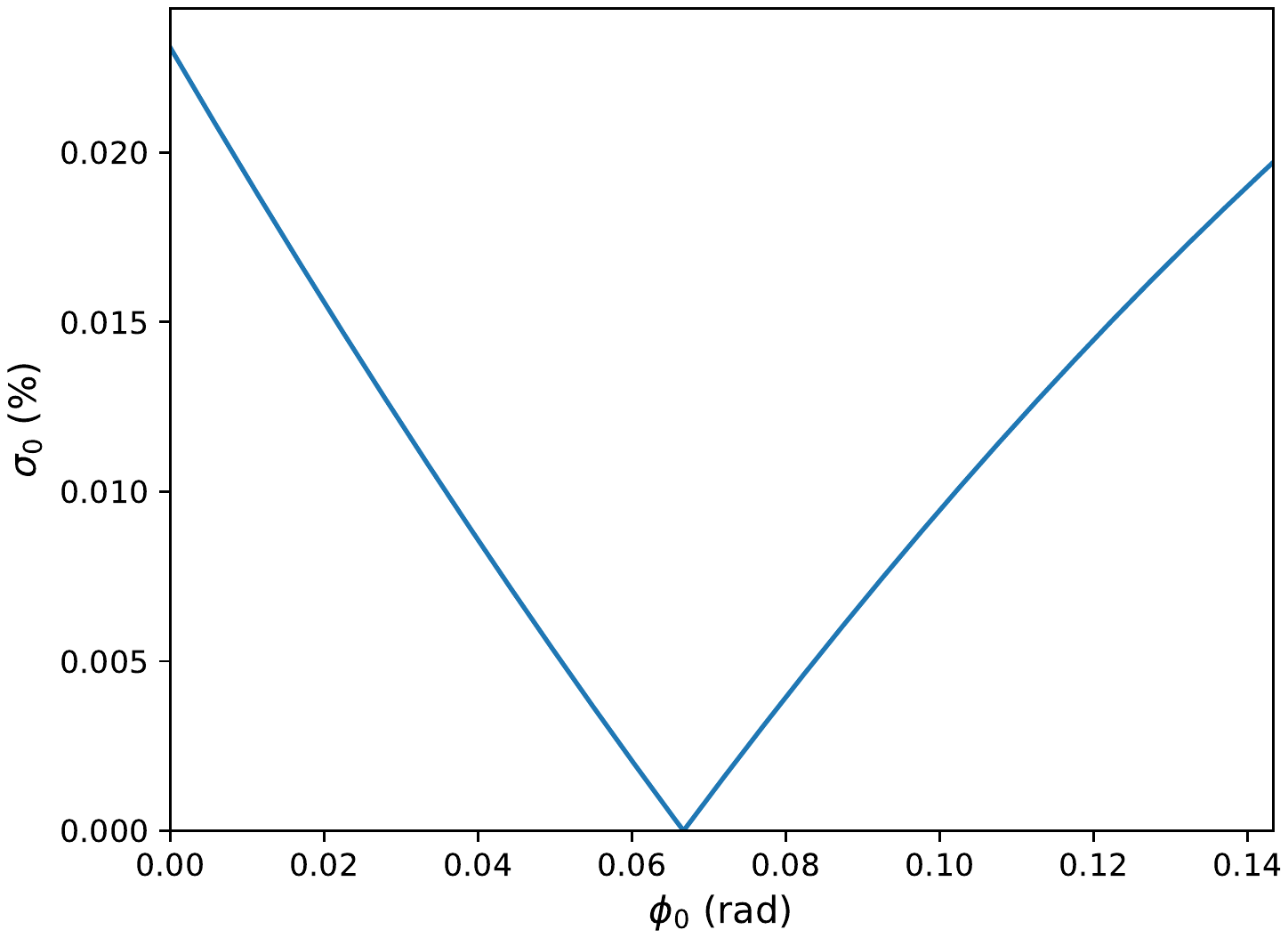}
			\caption{$r_A=3r$, $r_B=5r$, $2\theta=\frac{2\pi}{3}$}
		\end{subfigure}
	\caption{}
\end{figure}
\section{Conclusions.}
It was shown that it is possible in principle to determine the point of \emph{spherical specular reflection} by implementing a \emph{geometric algorithm} involving only \emph{ruler} and \emph{compass constructions}. Numerical results for the error involved in approximating the actual solution with the first iteration of the algorithm provided convincing evidence for it's use in applications.
\section*{Acknowledgments.}
The author wishes to thank K. Blekos for helpful discussions and remarks as well as G. Dasios for pointing out the problem of the spherical specular reflection point.
\appendix
\section{Solutions to the quartic equation.}
In Figure 9 we present the points on the surface of the sphere corresponding to the solutions of Equation \ref{eq1} where we have also included the original solution for completeness. Considerations similar to those in Section 2 will convince that in each case $(AMC)$ and $(BMD)$ are similar triangles and lead to the same equation. The designations \emph{internal-external} refer to the cases where the reflection angles with respect to the ray connecting the origin to $M$ lie within and outside the circle of radius $r$ respectively. Note how in cases (a) and (c) the angles at $O$ sum to $2\theta$ while those in cases (b) and (d) to $\pi-2\theta$.
\begin{figure}\label{fig8}
	\begin{subfigure}{0.4\textwidth}
		\includegraphics[scale=0.4]{ext-ext.JPG}
		\caption{External reflection.}
	\end{subfigure}
	\begin{subfigure}{0.4\textwidth}
		\includegraphics[scale=0.4]{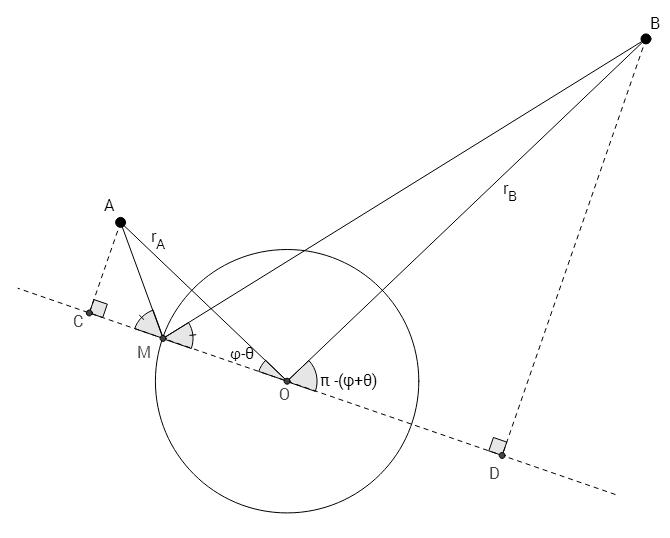}
		\caption{External-Internal reflection.}
	\end{subfigure}
	\begin{subfigure}{0.4\textwidth}
		\includegraphics[scale=0.4]{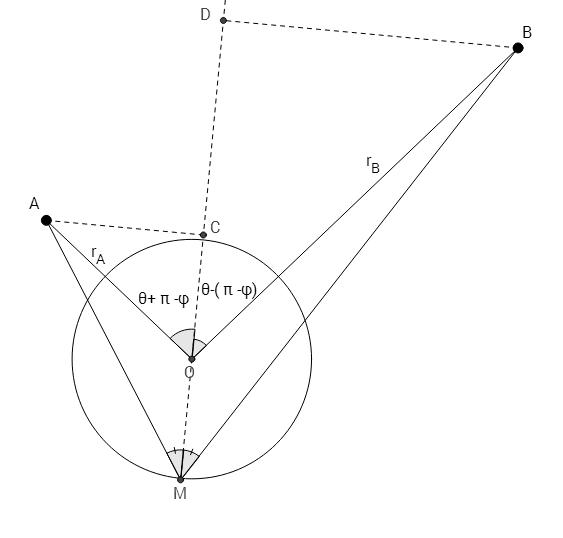}
		\caption{Internal reflection.}
	\end{subfigure}\hfill
	\begin{subfigure}{0.4\textwidth}
		\includegraphics[scale=0.4]{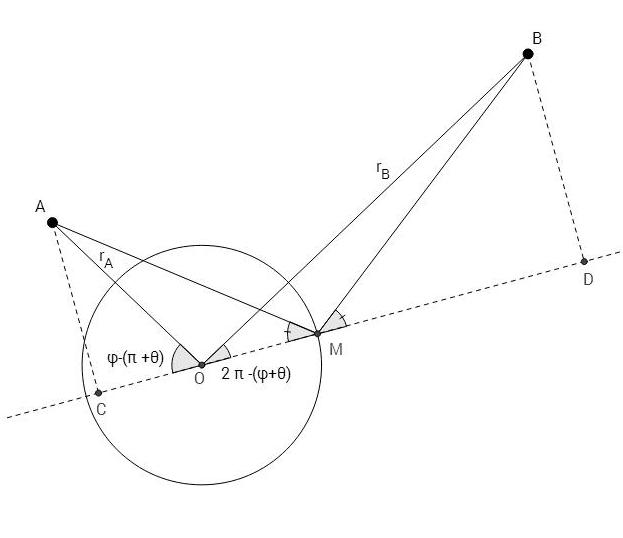}
		\caption{Internal-External reflection.}
	\end{subfigure}
	\caption{}
\end{figure}
\bibliographystyle{unsrt}
\bibliography{reflection}
\end{document}